\documentclass[11pt]{article}
\usepackage{amsfonts}
\title{THE SECOND PINCHING THEOREM FOR HYPERSURFACES WITH CONSTANT MEAN CURVATURE IN A
SPHERE \footnote{
\indent Research supported by the Chinese NSF, Grant
No. 11071211, 10771187; the Trans-Century Training Programme
Foundation for Talents by the Ministry of Education of China.}}
\author{HONG-WEI XU AND ZHI-YUAN XU}
\date{}
\begin{document}
\maketitle
\begin{abstract}  We generalize the second pinching theorem for minimal hypersurfaces in a sphere due to
Peng-Terng, Wei-Xu, Zhang, and Ding-Xin to the case of hypersurfaces
with small constant mean curvature. Let $M^n$ be a compact
hypersurface with constant mean curvature $H$ in $\mathbb{S}^{n+1}$.
Denote by $S$ the squared norm of the second fundamental form of
$M$. We prove that there exist two positive constants $\gamma(n)$
and $\delta(n)$ depending only on $n$ such that if
$|H|\leq\gamma(n)$ and $\beta(n,H)\leq S\leq\beta(n,H)+\delta(n)$,
then $S\equiv\beta(n,H)$ and $M$ is one of the following cases: (i)
$\mathbb{S}^{k}(\sqrt{\frac{k}{n}})\times
\mathbb{S}^{n-k}(\sqrt{\frac{n-k}{n}})$, $\,1\le k\le n-1$; (ii)
$\mathbb{S}^{1}(\frac{1}{\sqrt{1+\mu^2}})\times
\mathbb{S}^{n-1}(\frac{\mu}{\sqrt{1+\mu^2}})$. Here
$\beta(n,H)=n+\frac{n^3}{2(n-1)}H^2+\frac{n(n-2)}{2(n-1)}\sqrt{n^2H^4+4(n-1)H^2}$
and $\mu=\frac{n|H|+\sqrt{n^2H^2+4(n-1)}}{2}$.
\end{abstract}

\section*{1. Introduction}
\hspace*{5mm}Let $M^n$ be an $n$-dimensional compact hypersurface
with constant mean curvature $H$ in an $(n+1)$-dimensional unit
sphere $\mathbb{S}^{n+1}$. Denote by $S$ the squared length of the
second fundamental form of $M$ and $R$ its scalar curvature. Then
$R=n(n-1)+n^2H^{2}-S$. When $H=0$, the famous pinching theorem due
to Simons, Lawson, and Chern, do Carmo and Kobayashi (\cite{CDK},
\cite{L}, \cite{S}) says that if $S\leq n$, then $S\equiv0$ or
$S\equiv n$, i.e., $M$ must be the great sphere $\mathbb{S}^{n}$ or
the Clifford torus $\mathbb{S}^{k}(\sqrt{\frac{k}{n}})\times
\mathbb{S}^{n-k}(\sqrt{\frac{n-k}{n}})$, $\,1\le k\le n-1$. Further
discussions have been carried out by many other authors (see
\cite{D}, \cite{LL}, \cite{SCY}, \cite{X1}, \cite{X2}, \cite{Y},
etc.). In 1970's, Chern proposed the following conjectures.\\
\textbf{Chern Conjecture I.} \textit{ Let $M$ be a compact minimal
hypersurface with constant scalar curvature in $\mathbb{S}^{n+1}$.
Then the possible values form a discrete set. In particular, if
$n\le
S\le2n$, then $S=n$, or $S=2n.$}\\
\textbf{Chern Conjecture II.} \textit{ Let $M$ be a compact minimal
hypersurface in $\mathbb{S}^{n+1}$. If $n\le S\le2n$, then $S\equiv
n$, or $S\equiv 2n.$}
\par In 1983, Peng and Terng made breakthrough on the Chern
conjectures I and II. They \cite{PT1} proved that if $M$ is a
compact minimal hypersurface with constant scalar curvature in the
unit sphere $\mathbb{S}^{n+1}$, and if $n\leq S\leq
n+\frac{1}{12n}$, then $S=n$. Moreover, Peng and Terng \cite{PT2}
proved that if $M$ is a compact minimal hypersurface in the unit
sphere $\mathbb{S}^{n+1}$, and if $n\le 5$ and $n\leq S\leq
n+\tau_1(n)$, where $\tau_1(n)$ is a positive constant depending
only on $n$, then $S\equiv n$. During the past two decades, there
have been some important progress on these aspects(see \cite{C},
\cite{CI}, \cite{CY}, \cite{DX}, \cite{SY}, \cite{WX}, \cite{Z},
etc.).  In 1993, Chang \cite{C} solved Chern Conjecture I for the
case of dimension 3. In \cite{CI} and \cite{CY}, Cheng, Ishikawa and
Yang obtained some interesting results on the Chern conjectures.
\par In 2007, Suh-Yang and Wei-Xu made some progress on Chern Conjectures,
respectively. Suh and Yang \cite{SY} proved that if $M$ is a compact
minimal hypersurface with constant scalar curvature in
$\mathbb{S}^{n+1}$, and if $n\leq S\leq n+ \frac{3}{7}n$, then $S=n$
and $M$ is a minimal Clifford torus. Meanwhile, Wei and Xu \cite{WX}
proved that if $M$ is a compact minimal hypersurface in
$\mathbb{S}^{n+1}$, $n=6, 7$, and if $n\leq S\leq n+ \tau_2(n)$,
where $\tau_2(n)$ is a positive constant depending only on $n$, then
$S\equiv n$ and $M$ is a minimal Clifford torus. Later, Zhang
\cite{Z} extended the second pinching theorem due to Peng-Terng
\cite{PT2} and Wei-Xu \cite{WX} to $8$-dimensional compact minimal
hypersurfaces in a unit sphere. Recently Ding and Xin \cite{DX}
obtained the following pinching theorem for $n$-dimensional minimal
hypersurfaces in a
sphere.\\\\
\textbf{Theorem A}. \emph{Let $M$ be an $n$-dimensional compact
minimal hypersurface in a unit sphere $\mathbb{S}^{n+1}$, and $S$
the squared length of the second fundamental form of $M$. Then there
exists a positive constant $\tau(n)$ depending only on n such that
if $n\leq S\leq n+\tau(n)$, then $S\equiv n$, i.e., $M$ is a
Clifford torus.}\\\\
\hspace*{5mm}The pinching phenomenon for hypersurfaces of constant
mean curvature in spheres is much more complicated than the minimal
hypersurface case (see \cite{X1}, \cite{X3}). In \cite{X1}, Xu
proved the following pinching theorem for
submanifolds with parallel mean curvature in a sphere.\\\\
 \textbf{Theorem
B}. \emph{ Let $M$ be an $n$-dimensional compact submanifold with
parallel mean curvature vector $(H\neq0)$ in an $(n+p)$-dimensional
unit sphere $\mathbb{S}^{n+p}$. If $S\leq \alpha(n,H)$, then either
$M$ is pseudo-umbilical, or $S\equiv \alpha(n,H)$ and $M$ is the
isoparametric hypersurface
$\mathbb{S}^{n-1}\big(\frac{1}{\sqrt{1+\lambda^2}}\big)\times
\mathbb{S}^{1}\big(\frac{\lambda}{\sqrt{1+\lambda^2}}\big)$ in a
great sphere $\mathbb{S}^{n+1}$. In particular, if $M$ is a compact
hypersurface with constant mean curvature $H(\neq0)$ in
$\mathbb{S}^{n+1}$, then $M$ is either a totally umbilical sphere
$\mathbb{S}^n\big(\frac{1}{\sqrt{1+H^2}}\big)$, or a Clifford
hypersurface
$\mathbb{S}^{n-1}\big(\frac{1}{\sqrt{1+\lambda^2}}\big)\times
 \mathbb{S}^1\big(\frac{\lambda}{\sqrt{1+\lambda^2}}\big)$. Here
$\alpha(n,H)=n+\frac{n^3H^2}{2(n-1)}
-\frac{n(n-2)|H|}{2(n-1)}\sqrt{n^2H^2+4(n-1)}$
and $\lambda=\frac{n|H|+\sqrt{n^2H^2+4(n-1)}}{2(n-1)}.$}\\\\
\hspace*{5mm} In \cite{XT}, Xu and Tian generalized Suh-Yang's
pinching theorem \cite{SY} to the case where $M$ is a compact
hypersurface with constant scalar curvature and small constant mean
curvature in $\mathbb{S}^{n+1}$. The following second pinching
theorem for hypersurfaces with small constant mean curvature was
proved for $n\le 7$ by Cheng-He-Li \cite{CHL} and Xu-Zhao
\cite{XZ} respectively, and for $n=8$ by Xu \cite{XU}.\\\\
\textbf{Theorem C}. \emph{Let $M$ be an $n$-dimensional compact
hypersurface with constant mean curvature $H(\neq0)$ in a unit
sphere $\mathbb{S}^{n+1}$, $n\leq8$. There exist two positive
constants $\gamma_0(n)$ and $\delta_0(n)$ depending only on $n$ such
that if $|H|\leq\gamma_0(n)$, and $\beta(n,H)\leq
S<\beta(n,H)+\delta_0(n)$, then $S\equiv\beta(n,H)$ and
$M=\mathbb{S}^{1}(\frac{1}{\sqrt{1+\mu^2}})\times
\mathbb{S}^{n-1}(\frac{\mu}{\sqrt{1+\mu^2}})$. Here
$\beta(n,H)=n+\frac{n^3}{2(n-1)}H^2+\frac{n(n-2)}{2(n-1)}\sqrt{n^2H^4+4(n-1)H^2}$
and $\mu=\frac{n|H|+\sqrt{n^2H^2+4(n-1)}}{2}$.}\\\\
\hspace*{5mm} In this paper, we prove the second pinching theorem
for $n$-dimensional hypersurfaces with constant mean curvature,
which is a generalization of Theorems A and
C.\\\\
\textbf{Main Theorem}. \emph{Let $M$ be an $n$-dimensional compact
hypersurface with constant mean curvature $H$ in a unit sphere
$\mathbb{S}^{n+1}$. There exist two positive constants $\gamma(n)$
and $\delta(n)$ depending only on $n$ such that if
$|H|\leq\gamma(n)$, and $\beta(n,H)\leq S\leq\beta(n,H)+\delta(n)$,
then $S\equiv\beta(n,H)$ and $M$ is one of the following cases: (i)
$\mathbb{S}^{k}(\sqrt{\frac{k}{n}})\times
\mathbb{S}^{n-k}(\sqrt{\frac{n-k}{n}})$, $\,1\le k\le n-1$;  (ii)
$\mathbb{S}^{1}(\frac{1}{\sqrt{1+\mu^2}})\times
\mathbb{S}^{n-1}(\frac{\mu}{\sqrt{1+\mu^2}})$. Here
$\beta(n,H)=n+\frac{n^3}{2(n-1)}H^2+\frac{n(n-2)}{2(n-1)}\sqrt{n^2H^4+4(n-1)H^2}$
and $\mu=\frac{n|H|+\sqrt{n^2H^2+4(n-1)}}{2}$.}\\\\

\section*{2. Preliminaries}
\hspace*{5mm} Let $M^{n}$ be an $n$-dimensional compact hypersurface
with constant mean curvature in a unit sphere $\mathbb{S}^{n+1}$. We
shall make use of the following convention on the range of indices.
                     $$1 \leq A, B, C, \ldots, \leq n+1, \ \  1\leq i, j, k, \ldots, \leq n.$$
For an arbitrary fixed point $x\in M\subset \mathbb{S}^{n+1}$, we
choose an orthonormal local frame field $\{e_{A}\}$ in
$\mathbb{S}^{n+1}$ such that $e_{i}$'s are tangent to $M$. Let
$\{\omega_{A}\}$ be the dual frame fields of $\{e_{A}\}$ and
$\{\omega_{AB}\}$ the connection 1-forms of $\mathbb{S}^{n+1}$.
Restricting to $M$, we have
\begin{equation}
\omega_{n+1 i}=\sum\limits_{j}h_{ij}\omega_{j}, \, h_{ij}=h_{ji}.
\end{equation}
\hspace*{5mm} Let $h$ be the second fundamental form of $M$. Denote
by $R$, $H$ and $S$ the scalar curvature, mean curvature and squared
length of the second fundamental form of $M$, respectively. Then we
have
\begin{equation}
h=\sum\limits_{i,j}h_{ij}\omega_{i}\otimes\omega_{j},
\end{equation}
\begin{equation}
S=\sum\limits_{i,j}h_{ij}^{2}, \,\,
H=\frac{1}{n}\sum\limits_{i}h_{ii},
\end{equation}
\begin{equation}
R=n(n-1)+n^2H^{2}-S.
\end{equation}
We choose $e_{n+1}$ such that
$H=\frac{1}{n}\sum\limits_{i}h_{ii}\ge0.$ Denote by $h_{ijk}$,
$h_{ijkl}$ and $h_{ijklm}$ the first, second and third covariant
derivatives of the second fundamental tensor $h_{ij}$, respectively.
Then we have
\begin{equation}
\nabla
h=\sum\limits_{i,j,k}h_{ijk}\omega_{i}\otimes\omega_{j}\otimes\omega_{k},\,\,
h_{ijk}=h_{ikj},
\end{equation}
\begin{equation}
h_{ijkl}=h_{ijlk}+\sum\limits_{m}h_{mj}R_{mikl}+\sum\limits_{m}h_{im}R_{mjkl},
\end{equation}
\begin{equation}
h_{ijklm}=h_{ijkml}+\sum\limits_{r}h_{rjk}R_{rilm}+\sum\limits_{r}h_{irk}R_{rjlm}+\sum\limits_{r}h_{ijr}R_{rklm}.
\end{equation}
\hspace*{5mm}At each fixed point $x\in M$, we take orthonormal
frames $\{e_{i}\}$ such that $h_{ij}=\lambda_{i}\delta_{ij}$ for all
$i$, $j$. Then $\sum\limits_{i}\lambda_{i}=nH$ and
$\sum\limits_{i}\lambda_{i}^{2}=S$. By a direct computation, we have
\begin{equation}
\frac{1}{2}\Delta S=S(n-S)-n^2H^2+nHf_3+|\nabla h|^2,
\end{equation}
\begin{eqnarray}
\frac{1}{2}\Delta|\nabla h|^2&=&\nonumber (2n+3-S)|\nabla
h|^2-\frac{3}{2}|\nabla S|^2+|\nabla^2 h|^2\\&
&\nonumber+\sum\limits_{i,j,k,l,m}(6h_{ijk}h_{ilm}h_{jl}h_{km}-3h_{ijk}h_{ijl}h_{km}h_{ml})
+3nH\sum\limits_{i,j,k,l}h_{ijk}h_{jlk}h_{li}\\&=& (2n+3-S)|\nabla
h|^2-\frac{3}{2}|\nabla S|^2+|\nabla^2 h|^2+3(2B-A)+3nHC,
\end{eqnarray}
where $$f_k=\sum\limits_{i}\lambda_{i}^k,\,\,
A=\sum\limits_{i,j,k}h_{ijk}^2\lambda_{i}^2,\,\,
B=\sum\limits_{i,j,k}h_{ijk}^2\lambda_{i}\lambda_{j},\,\,
C=\sum\limits_{i,j,k}h_{ijk}^2\lambda_{i}.$$\\
Using a similar method as in \cite{PT1}, we obtain
\begin{equation}
h_{ijij}=h_{jiji}+t_{ij},
\end{equation}
\begin{equation}
|\nabla^2 h|^2\geq\frac{3}{4}\sum\limits_{i\neq
j}t_{ij}^2=\frac{3}{4}\sum\limits_{i,j}t_{ij}^2,
\end{equation}
and
\begin{equation}
3(A-2B)\leq aS|\nabla h|^2,
\end{equation}
where $t_{ij}=(\lambda_i-\lambda_j)(1+\lambda_i\lambda_j)$ and
$a=\frac{\sqrt{17}+1}{2}$. From (11), we have
\begin{equation}
|\nabla^2
h|^2\geq\frac{3}{2}[Sf_4-f_{3}^2-S^2-S(S-n)-n^2H^2+2nHf_3].
\end{equation}
By a computation, we obtain
\begin{eqnarray}
\frac{1}{3}\sum\limits_{i,j}h_{ij}(f_3)_{ij}&=&\nonumber\frac{1}{3}\sum\limits_{k}\lambda_k(f_3)_{kk}\\
&=&\nonumber\sum\limits_{k}\lambda_k(\sum\limits_{i}h_{iikk}\lambda_{i}^2+2\sum\limits_{i,j}h_{ijk}^2\lambda_i)\\
&=&\nonumber\sum\limits_{i,k}h_{iikk}\lambda_k\lambda_{i}^2+2\sum\limits_{i,j,k}h_{ijk}^2\lambda_i\lambda_k\\
&=&\nonumber\sum\limits_{i,k}[h_{kkii}+(\lambda_i-\lambda_k)(1+\lambda_i\lambda_k)]\lambda_k\lambda_{i}^2+2B\\
&=&\nonumber\sum\limits_{i}(\frac{S_{ii}}{2}-\sum\limits_{j,k}h_{ijk}^2)\lambda_{i}^2
+\sum\limits_{i,k}\lambda_{i}^2\lambda_{k}(\lambda_i-\lambda_k)(1+\lambda_i\lambda_k)+2B\\
&=&\sum\limits_{i,j,k}\frac{h_{ik}h_{kj}}{2}S_{ij}+nHf_3-S^2-f_{3}^2+Sf_4-(A-2B).
\end{eqnarray}
Since $\int_{M}\sum\limits_{i,j}h_{ij}(f_3)_{ij}dM=0$, we drive the
following integral formula.
\begin{eqnarray}
\int_{M}(A-2B)dM&=&\nonumber\int_{M}(nHf_3-S^2-f_{3}^2+Sf_4+\sum\limits_{i,j,k}\frac{h_{ik}h_{kj}}{2}S_{ij})dM\\
&=&\nonumber\int_{M}(nHf_3-S^2-f_{3}^2+Sf_4-\sum\limits_{i,j,k}(h_{ik}h_{kj})_j\frac{S_i}{2})dM\\
&=&\nonumber\int_{M}(nHf_3-S^2-f_{3}^2+Sf_4-\sum\limits_{i,j,k}h_{ikj}h_{kj}\frac{S_i}{2}-\sum\limits_{i,j,k}h_{ik}h_{kjj}\frac{S_i}{2})dM\\
&=&\nonumber\int_{M}(nHf_3-S^2-f_{3}^2+Sf_4-\sum\limits_{i,j,k}h_{ikj}h_{kj}\frac{S_i}{2})dM\\
&=&\int_{M}(nHf_3-S^2-f_{3}^2+Sf_4-\frac{|\nabla S|^2}{4})dM.
\end{eqnarray}

\section*{3. Proof of Main Theorem}

\hspace*{5mm} The key to the proof of Main Theorem is to establish
some integral equalities and inequalities on the second
fundamental form of $M$ and its covariant derivatives by the parameter method.\\
\hspace*{5mm} To simplify the computation, we introduce the
tracefree second fundamental form
$\phi=\sum\limits_{i,j}\phi_{ij}\omega_{i}\otimes\omega_{j}$, where
$\phi_{ij}=h_{ij}-H\delta_{ij}$. If $h_{ij}=\lambda_i\delta_{ij}$,
then $\phi_{ij}=\mu_i\delta_{ij}$, where $\mu_i=\lambda_i-H$.
Putting $\Phi=|\phi|^2$ and $\bar{f}_k=\sum\limits_{i}\mu_i^k$, we
get $\Phi=S-nH^2$, $f_3=\bar{f}_3+3H\Phi+nH^3$ and
$f_4=\bar{f}_4+4H\bar{f}_3+6H^2\Phi+nH^4$. From (8), we obtain
\begin{eqnarray}
\frac{1}{2}\Delta\Phi&=&\nonumber S(n-S)-n^2H^2+nHf_3+|\nabla h|^2\\
&=&\nonumber-\Phi^2+n\Phi+nH\bar{f}_3+nH^2\Phi+|\nabla\phi|^2\\
&=&-F(\Phi)+|\nabla\phi|^2,
\end{eqnarray}
where $F(\Phi)=\Phi^2-n\Phi-nH^2\Phi-nH\bar{f}_3.$ Therefore, we
have
\begin{equation}
|\nabla\Phi|^2=\frac{1}{2}\Delta(\Phi)^2-\Phi\Delta\Phi=\frac{1}{2}\Delta(\Phi)^2+2\Phi
F(\Phi)-2\Phi|\nabla\phi|^2,
\end{equation}
and
\begin{equation}
\int_MF(\Phi)dM=\int_M|\nabla\phi|^2dM.
\end{equation}
\textbf{Lemma 1}.(See \cite{X1}) \emph{Let $a_1$, $a_2$, ..., $a_n$
be real numbers satisfying $\sum\limits_{i}a_i=0$ and
$\sum\limits_{i}a_i^2=a$.
Then$$|\sum\limits_{i}a_i^3|\leq\frac{n-2}{\sqrt{n(n-1)}}a^{\frac{3}{2}},$$
and the equality holds if and only if at least $n-1$ numbers of
 $a_i$'s are same with each other.}\\\\
\hspace*{5mm} From Lemma 1, we get
\begin{eqnarray}
F(\Phi)&\geq&\nonumber\Phi^2-n\Phi-nH^2\Phi-\frac{n(n-2)H\Phi^{\frac{3}{2}}}{\sqrt{n(n-1)}}\\
&=&\nonumber\Phi\Big[\Phi-\frac{n(n-2)H\Phi^{\frac{1}{2}}}{\sqrt{n(n-1)}}-n(1+H^2)\Big]\\
&\geq&0,
\end{eqnarray}
provided
$$\Phi\geq\beta_0(n,H):=n+\frac{n^3}{2(n-1)}H^2+\frac{n(n-2)}{2(n-1)}\sqrt{n^2H^4+4(n-1)H^2}-nH^2.$$
Moreover, $F(\Phi)=0$ if and only if $\Phi=\beta_0(n,H)$.\\
\hspace*{5mm}Set$$G=\sum\limits_{i,j}(\lambda_i-\lambda_j)^2(1+\lambda_i\lambda_j)^2.$$
Then we have
\begin{equation}
G=2[Sf_4-f_3^2-S^2-S(S-n)+2nHf_3-n^2H^2].
\end{equation}
This together with (8) and (15) implies
\begin{equation}
\frac{1}{2}\int_MGdM=\int_M[(A-2B)-|\nabla h|^2+\frac{1}{4}|\nabla
S|^2]dM.
\end{equation}
\textbf{Lemma 2}. \emph{Let $M$ be an $n(\geq4)$-dimensional compact
hypersurface with constant mean curvature in $\mathbb{S}^{n+1}$. If
$S\geq\beta(n,H)$, then we have
$$3(A-2B)\leq 2S|\nabla h|^2+C_1(n)|\nabla h|^2G^{\frac{1}{3}},$$
where
$C_1(n)=(\sqrt{17}-3)[6(\sqrt{17}+1)]^{-\frac{1}{3}}(\frac{2}{\sqrt{17}}-
\frac{\sqrt{2}}{17}-\frac{1}{n})^{-\frac{2}{3}}$.
}\\
\noindent\textbf{Proof.} We derive the estimate above at each fixed
point $x\in M$. If $\lambda_j^2-4\lambda_i\lambda_j\leq2S$ for all
$i\neq j$, then we get the desired estimate immediately. Otherwise,
we assume that there exist $i\neq j$, such that
$\lambda_j^2-4\lambda_i\lambda_j=tS>2S$.\\
We get
\begin{equation}
S\geq\lambda_i^2+\lambda_j^2=(\frac{tS-\lambda_j^2}{4\lambda_j})^2+\lambda_j^2.
\end{equation}
Then
\begin{equation}
\lambda_j^2\leq\frac{1}{17}(t+8+4\sqrt{4+t-t^2})S,\,\,\,\,
2<t\leq\frac{\sqrt{17}+1}{2},
\end{equation}
which implies
\begin{equation}
-\lambda_i\lambda_j\geq\frac{1}{17}(4t-2-\sqrt{4+t-t^2})S\geq0.26S>\frac{S}{n}\geq1.
\end{equation}
On the other hand, we have
\begin{equation}
(\lambda_i-\lambda_j)^2=(\frac{\lambda_j}{2}+\lambda_i)^2+\frac{3}{4}(\lambda_j^2-4\lambda_i\lambda_j)\geq\frac{3t}{4}S.
\end{equation}
By the definition of $G$, we get
\begin{eqnarray}
G&\geq&\nonumber2(\lambda_i-\lambda_j)^2(1+\lambda_i\lambda_j)^2\\
&\geq&\nonumber\frac{3t}{2}S(1+\lambda_i\lambda_j)^2\\
&\geq&\nonumber\frac{3t}{2}S(-\lambda_i\lambda_j-\frac{S}{n})^2\\
&\geq&\frac{3t}{2}\Big[\frac{1}{17}(4t-2-\sqrt{4+t-t^2})-\frac{1}{n}\Big]^2S^3.
\end{eqnarray}
We define an auxiliary function
$$\zeta(t)=\frac{t}{(t-2)^3}\Big[\frac{1}{17}(4t-2-\sqrt{4+t-t^2})-\frac{1}{n}\Big]^2,\,\,\,\, 2<t\leq\frac{\sqrt{17}+1}{2}.$$
Then we have
\begin{eqnarray}
\zeta(t)&\geq&\nonumber\frac{t}{(t-2)^3}\Big[\frac{1}{17}(4t-2-\sqrt{2})-\frac{1}{n}\Big]^2\\
&\geq&\nonumber
\inf_{2<t\leq\frac{\sqrt{17}+1}{2}}\frac{t}{(t-2)^3}\Big[\frac{1}{17}(4t-2-\sqrt{2})-\frac{1}{n}\Big]^2\\
&=&\frac{4(\sqrt{17}+1)}{(\sqrt{17}-3)^3}\Big(\frac{2}{\sqrt{17}}-\frac{\sqrt{2}}{17}-\frac{1}{n}\Big)^2.
\end{eqnarray}
Hence
\begin{eqnarray}
(\lambda_j^2-4\lambda_i\lambda_j-2S)^3&=&\nonumber(t-2)^3S^3\\
&\leq&\nonumber\frac{2G}{3\zeta(t)}\\
&\leq&\nonumber\frac{(\sqrt{17}-3)^3}{6(\sqrt{17}+1)}(\frac{2}{\sqrt{17}}-\frac{\sqrt{2}}{17}-\frac{1}{n})^{-2}G\\
&=&(C_1(n)G^{\frac{1}{3}})^3.
\end{eqnarray}
This implies
\begin{eqnarray}
3(A-2B)&\leq&\nonumber\sum\limits_{i,j,k
\,\,distinct}[2(\lambda_i^2+\lambda_j^2+\lambda_k^2)-(\lambda_i+\lambda_j+\lambda_k)^2]h_{ijk}^2+3\sum\limits_{i\neq
j}(\lambda_j^2-4\lambda_i\lambda_j)h_{iij}^2\\
&\leq&\nonumber2S\sum\limits_{i,j,k
\,\,distinct}h_{ijk}^2+3\sum\limits_{i\neq
j}h_{iij}^2(2S+C_1(n)G^\frac{1}{3})\\
&\leq&2S|\nabla h|^2+C_1(n)|\nabla h|^2G^{\frac{1}{3}}.
\end{eqnarray}
\\\\
\noindent\textbf{Proof of Main Theorem}.(i) When $H=0$, the
assertion follows from Theorem A. \\
(ii) When $H\neq0$, the assertion for lower dimensional cases $(n\le
8)$ was verified in \cite{CHL}, \cite{XZ} and \cite{XU}. We consider
the case for $n\ge 4.$ From (10) and (11), we see that
$G=\sum\limits_{i,j}t_{ij}^2$ and $|\nabla^2 h|^2\geq\frac{3}{4}G$.
Let $0<\theta<1$, we have
\begin{equation}
\int_M|\nabla^2
h|^2dM\geq\Big[\frac{3(1-\theta)}{4}+\frac{3\theta}{4}\Big]\int_MGdM.
\end{equation}
From (9), (21), Lemma 2 and Young's inequality, we drive the
following inequality.
\begin{eqnarray}
\frac{3(1-\theta)}{4}\int_MGdM&\leq&\nonumber\int_M\Big[(S-2n-3)|\nabla
h|^2+\frac{3}{2}|\nabla S|^2+3(A-2B)-3nHC-\frac{3\theta}{4}G\Big]dM\\
&=&\nonumber\int_M(S-2n-3+\frac{3\theta}{2})|\nabla
h|^2dM+(3-\frac{3\theta}{2})\int_M(A-2B)dM\\
& &\nonumber+(\frac{3}{2}-\frac{3\theta}{8})\int_M|\nabla
S|^2dM-3nH\int_MCdM\\
&\leq&\nonumber\int_M(S-2n-3+\frac{3\theta}{2})|\nabla
h|^2dM+(1-\frac{\theta}{2})\int_M(2S|\nabla h|^2\\
& &\nonumber+C_1(n)|\nabla
h|^2G^{\frac{1}{3}})dM+(\frac{3}{2}-\frac{3\theta}{8})\int_M|\nabla
S|^2dM-3nH\int_MCdM\\
&\leq&\nonumber\int_M\Big[(3-\theta)S-2n-3+\frac{3\theta}{2}\Big]|\nabla
h|^2dM+\frac{3(1-\theta)}{4}\int_MGdM\\
& &\nonumber+C_2(n,\theta)\int_M|\nabla
h|^3dM+(\frac{3}{2}-\frac{3\theta}{8})\int_M|\nabla
S|^2dM\\
& &-3nH\int_MCdM,
\end{eqnarray}
where
$C_2(n,\theta)=\frac{4}{9}C_1(n)^{\frac{3}{2}}(1-\frac{\theta}{2})^{\frac{3}{2}}(1-\theta)^{-\frac{1}{2}}$.\\\\
\hspace*{5mm}Let $\epsilon>0$, from (16), we get
\begin{eqnarray}
\int_M|\nabla h|^3dM&=&\nonumber\int_M|\nabla\phi|^3dM\\
&=&\nonumber\int_M|\nabla\phi|(F(\Phi)+\frac{1}{2}\Delta\Phi)dM\\
&=&\nonumber\int_MF(\Phi)|\nabla\phi|dM-\frac{1}{2}\int_M\nabla|\nabla\phi|\cdot\nabla\Phi dM\\
&\leq&\int_MF(\Phi)|\nabla\phi|dM+\epsilon\int_M|\nabla^2\phi|^2dM+\frac{1}{16\epsilon}\int_M|\nabla\Phi|^2dM.
\end{eqnarray}
Since
\begin{equation}
|C|\leq \sqrt{S}|\nabla h|^2,
\end{equation}
we have
\begin{eqnarray}
0&\leq&\nonumber\int_M[(3+3\sqrt{n}H-\theta)(\Phi+nH^2)-2n-3+\frac{3\theta}{2}]|\nabla\phi|^2dM\\
&
&\nonumber+C_2(n,\theta)[\int_MF(\Phi)|\nabla\phi|dM+\epsilon\int_M|\nabla^2\phi|^2dM
+\frac{1}{16\epsilon}\int_M|\nabla\Phi|^2dM]\\
& &+(\frac{3}{2}-\frac{3\theta}{8})\int_M|\nabla\Phi|^2dM.
\end{eqnarray}
Substituting (12) and (33) into (9), we have
\begin{eqnarray}
\int_M|\nabla^2\phi|^2dM&=&\nonumber\int_M|\nabla^2h|^2dM\\
&\leq&\nonumber\int_M[(S-2n-3)|\nabla h|^2+\frac{3}{2}|\nabla
S|^2+aS|\nabla h|^2-3nHC]dM\\
&\leq&\int_M[(a+1+3\sqrt{n}H)S-2n-3]|\nabla\phi|^2dM+\frac{3}{2}\int_M|\nabla
S|^2dM.
\end{eqnarray}
Combining (16) and (17), we have
\begin{eqnarray}
\int_M\frac{1}{2}|\nabla\Phi|^2dM&=&\nonumber\int_M\Phi
F(\Phi)dM-\int_M\Phi|\nabla\phi|^2dM+\beta_0(n,H)\int_M|\nabla\phi|^2dM\\
& &\nonumber-\beta_0(n,H)\int_MF(\Phi)dM\\
&=&\int_M(\Phi-\beta_0(n,H))F(\Phi)dM+\int_M(\beta_0(n,H)-\Phi)|\nabla\phi|^2dM.
\end{eqnarray}
Hence
\begin{eqnarray}
0&\leq&\nonumber\int_M\Big\{\Big[3+3\sqrt{n}H-\theta+\epsilon C_2(n,\theta)(a+1+3\sqrt{n}H)\Big](\Phi-\beta_0(n,H))\\
& &\nonumber+\beta(n,H)\Big[3+3\sqrt{n}H-\theta+\epsilon C_2(n,\theta)(a+1+3\sqrt{n}H)\Big]\\
&
&\nonumber-2\Big(\frac{3}{2}-\frac{3\theta}{8}+\frac{C_2(n,\theta)}{16\epsilon}+
\frac{3\epsilon
C_2(n,\theta)}{2}\Big)(\Phi-\beta_0(n,H))\\
& &\nonumber-2n-3+\frac{3\theta}{2}-\epsilon
C_2(n,\theta)(2n+3)\Big\}|\nabla\phi|^2dM\\
&
&\nonumber+2\Big(\frac{3}{2}-\frac{3\theta}{8}+\frac{C_2(n,\theta)}{16\epsilon}+
\frac{3\epsilon
C_2(n,\theta)}{2}\Big)\int_M(\Phi-\beta_0(n,H))F(\Phi)dM\\
& &\nonumber+C_2(n,\theta)\int_MF(\Phi)|\nabla\phi|dM\\
&=&\nonumber\int_M\Big\{D(n,H)\Big[3+3\sqrt{n}H-\theta+\epsilon
C_2(n,\theta)(a+1+3\sqrt{n}H)\Big]\\
&
&\nonumber+(1-\theta)n-3+\frac{3\theta}{2}+3n^{\frac{3}{2}}H+\epsilon
C_2(n,\theta)(a
n+3n^{\frac{3}{2}}H-n-3)\Big\}|\nabla\phi|^2dM\\
&
&\nonumber-\Big(\frac{\theta}{4}+\frac{C_2(n,\theta)}{8\epsilon}-3\sqrt{n}H+\epsilon
C_2(n,\theta)(2-a-3\sqrt{n}H)\Big)\int_M(\Phi-\beta_0(n,H))|\nabla\phi|^2dM\\
&
&\nonumber+\Big(3-\frac{3\theta}{4}+\frac{C_2(n,\theta)}{8\epsilon}+3\epsilon
C_2(n,\theta)\Big)\int_M(\Phi-\beta_0(n,H))F(\Phi)dM\\
& &+C_2(n,\theta)\int_MF(\Phi)|\nabla\phi|dM,
\end{eqnarray}
where $\beta(n,H)=\beta_0(n,H)+nH^2$ and
$D(n,H)=\beta(n,H)-n$.\\\\
\hspace*{5mm}Note that
\begin{equation}
\frac{\theta}{4}+\frac{C_2(n,\theta)}{8\epsilon}-3\sqrt{n}H+\epsilon
C_2(n,\theta)(2-a-3\sqrt{n}H)\geq0,
\end{equation}
for all $\epsilon \in (0, \epsilon_1]$, where $\epsilon_1$ is some
positive constant. When $\beta(n,H)\leq S\leq\beta(n,H)+\epsilon^2$,
we obtain
\begin{eqnarray}
0&\leq&\nonumber\int_M[(1-\theta)n-3+\frac{3\theta}{2}+3n^{\frac{3}{2}}H+D(n,H)(3+3\sqrt{n}H-\theta)+O(\epsilon,\theta,H)]|\nabla\phi|^2dM\\
& &+C_2(n,\theta)\int_MF(\Phi)|\nabla\phi|dM,
\end{eqnarray}
where
\begin{eqnarray}
O(\epsilon,\theta,H)&=&\nonumber\epsilon
D(n,H)C_2(n,\theta)(a+1+3\sqrt{n}H)+\epsilon C_2(n,\theta)(a
n+3n^{\frac{3}{2}}H-n-3)\\
&
&\nonumber+\epsilon^2(3-\frac{3\theta}{4}+\frac{C_2(n,\theta)}{8\epsilon}+3\epsilon
C_2(n,\theta)).
\end{eqnarray}
On the other hand, we have
\begin{equation}
C_2(n,\theta)\int_MF(\Phi)|\nabla\phi|dM\leq\frac{3}{8}\int_MF(\Phi)dM+\frac{2C_2(n,\theta)^2}{3}\int_MF(\Phi)|\nabla\phi|^2dM.
\end{equation}
Using Lemma 1, we drive an upper bound for $F(\Phi)$.
\begin{eqnarray}
F(\Phi)&\leq&\nonumber\Phi^2-n\Phi-nH^2\Phi+\frac{n(n-2)H\Phi^{\frac{3}{2}}}{\sqrt{n(n-1)}}\\
&=&\nonumber\Phi\Big[\Phi+\frac{n(n-2)H\Phi^{\frac{1}{2}}}{\sqrt{n(n-1)}}-n(1+H^2)\Big]\\
&=&\frac{\Phi(\Phi^{\frac{1}{2}}+\beta_0(n,H)^{\frac{1}{2}})(\Phi-\alpha_0(n,H))}{\Phi^{\frac{1}{2}}+\alpha_0(n,H)^{\frac{1}{2}}},
\end{eqnarray}
where $\alpha_0(n,H)=\Big[\frac{-n(n-2)H+n\sqrt{n^2H^2+4n-4}}{2\sqrt{n(n-1)}}\Big]^2$.\\
\hspace*{5mm}When $\delta(n)\leq\epsilon^2$ and $\epsilon\leq1$, we
choose positive constant $\gamma_1(n)$ such that $n\leq\Phi\leq2n$
and $x_1\leq2\sqrt{n}$ for all $H\leq\gamma_1(n)$. We obtain
\begin{equation}
F(\Phi)\leq8n(\Phi-\alpha_0(n,H))\leq8n\Big(\epsilon^2+\frac{n(n-2)}{(n-1)}
\sqrt{n^2H^4+4(n-1)H^2}\Big).
\end{equation}
\hspace*{5mm}Let $\theta=\theta(n)=1-\frac{1}{8n}$. We choose
positive constants $\gamma_2(n)$ and $\gamma_3(n)$ such that
$3n^{\frac{3}{2}}H+D(n,H)(3+3\sqrt{n}H)\leq\frac{1}{8}$ for all
$H\leq\gamma_2(n)$, and $\frac{16n^2(n-2)}{(n-1)}
\sqrt{n^2\gamma_3(n)^4+4(n-1)\gamma_3(n)^2}\leq\frac{9}{16C_2(n,\theta(n))^2}$.\\
\hspace*{5mm}Take $\epsilon_2(n)=\Big[\frac{n(n-2)}{(n-1)}
\sqrt{n^2\gamma_3(n)^4+4(n-1)\gamma_3(n)^2}\Big]^{\frac{1}{2}}>0$.
Combining (39), (40) and (42), we obtain
\begin{equation}
\int_M[-\frac{1}{2}+O(\epsilon,\theta(n),H)]|\nabla\phi|^2dM\geq0,
\end{equation}
for all $H\leq\gamma(n)=\min\{\gamma_1(n), \gamma_2(n), \gamma_3(n)\}$ and $\epsilon\leq\min\{\epsilon_1, \epsilon_2(n)\}$.\\
\hspace*{5mm}For $\epsilon\leq1$, we have
\begin{eqnarray}
O(\epsilon,\theta(n),H)&\leq&\nonumber\epsilon
D(n,\gamma(n))C_2(n,\theta(n))(a+1+3\sqrt{n}\gamma(n))\\
& &\nonumber+\epsilon C_2(n,\theta(n))(a
n+3n^{\frac{3}{2}}\gamma(n))\\
&
&\nonumber+\epsilon(3-\frac{3\theta(n)}{4}+\frac{C_2(n,\theta(n))}{8}+3
C_2(n,\theta(n)))\\
&:=&\epsilon \eta(n),
\end{eqnarray}
where $a=\frac{\sqrt{17}+1}{2}$.\\
\hspace*{5mm}For $\epsilon\leq\epsilon_1(n)$, where
$\epsilon_1(n)=\frac{C_2(n,\theta(n))}{8[3\sqrt{n}\gamma(n)+
C_2(n,\theta(n))(a+3\sqrt{n}\gamma(n)-2)]}>0$,
$a=\frac{\sqrt{17}+1}{2}$, we have
\begin{equation}
\frac{C_2(n,\theta(n))}{8\epsilon}\geq3\sqrt{n}\gamma(n)+
C_2(n,\theta(n))(a+3\sqrt{n}\gamma(n)-2)-\frac{\theta(n)}{4}.
\end{equation}
So
$$\frac{\theta(n)}{4}+\frac{C_2(n,\theta(n))}{8\epsilon}-3\sqrt{n}H+\epsilon
C_2(n,\theta(n))(2-a-3\sqrt{n}H)\geq0.$$ Taking
$\delta(n)=\epsilon(n)^2$, where
$\epsilon(n)=\min\{1,\epsilon_1(n),\epsilon_2(n), \epsilon_3(n)\}$
and $\epsilon_3(n)=\frac{1}{3\eta(n)}$, we have $\delta(n)>0$. From
(43) and the assumption that $\beta(n,H)\leq
S\leq\beta(n,H)+\delta(n)$, we obtain $\nabla\phi=0$. This implies
$F(\Phi)=0$ and $\Phi=\beta_0(n,H)$.\\
\hspace*{5mm}By Lemma 1, we have
$$\lambda_1=\dots=\lambda_{n-1}=H-\sqrt{\frac{\beta(n,H)-nH^2}{n(n-1)}},$$
$$\lambda_n=H+\sqrt{\frac{(n-1)(\beta(n,H)-nH^2)}{n}}.$$
Therefore $M$ is the Clifford hypersurface
$$\mathbb{S}^{1}(\frac{1}{\sqrt{1+\mu^2}})\times \mathbb{S}^{n-1}(\frac{\mu}{\sqrt{1+\mu^2}})$$
in $\mathbb{S}^{n+1}$, where
$\mu=\frac{nH+\sqrt{n^2H^2+4(n-1)}}{2}$. This completes the proof of
Main
Theorem.\\\\
\hspace*{5mm}Finally we would like to propose the following problems.\\\\
\textbf{Open Problem A}. \emph{Let $M$ be an $n$-dimensional compact
hypersurface with constant mean curvature $H$ in the unit sphere
$\mathbb{S}^{n+1}$. Does there exist a positive constant $\delta(n)$
depending only on $n$ such that if $\beta(n,H)\leq
S\leq\beta(n,H)+\delta(n)$, then
$S\equiv\beta(n,H)$?}\\\\
\textbf{Open Problem B}. \emph{ For an $n$-dimensional compact
hypersurface $M^n$ with constant mean curvature $H$ in
$\mathbb{S}^{n+1}$, set
$\mu_k=\frac{n|H|+\sqrt{n^2H^2+4k(n-k)}}{2k}$. Suppose that
$\alpha(n,H)\leq S\leq \beta(n,H)$. Is it possible to prove that $M$
must be the isoparametric hypersurface
$S^k\big(\frac{1}{\sqrt{1+\mu_k^2}}\big)\times
S^{n-k}\big(\frac{\mu_k}{\sqrt{1+\mu_k^2}}\big)$,
$k=1,2,\cdots,n-1$?}\\\\
\hspace*{5mm} When $H=0$, the rigidity theorem due to Lawson
\cite{L}, Chern, do Carmo and Kobayashi \cite{CDK} provides an
affirmative answer for Open
Problem B.\\\\
\textbf{\large{Acknowledgement.}}\ We would like to thank Dr. En-Tao
Zhao for his helpful discussions. Thanks also to Professor Y. L. Xin
for sending us the reference \cite{DX}.

Center of Mathematical Sciences\

Zhejiang University\

Hangzhou 310027\

China\\\

e-mail address: xuhw@cms.zju.edu.cn;  srxwing@zju.edu.cn

\end{document}